\documentclass[11pt,a4paper]{article}%

\usepackage{amsmath}%
\usepackage{graphicx}%
\usepackage[sort]{natbib}%
\usepackage{epstopdf}
\usepackage[top=1.7in, bottom=1.7in, left=0.9in, right=0.9in]{geometry}

\newcommand{\keywords}[1]{\par\addvspace\baselineskip\noindent{\em Some key words}:\enspace\ignorespaces#1}

\title{On the Mathematics of the Jeffreys--Lindley Paradox}
\date{}
\author{Cristiano Villa and Stephen Walker\\ School of Mathematics, Statistics and Actuarial Science - University of Kent\\Division of Statistics and Scientific Computation - University of Texas at Austin}
\def\half{\hbox{$1\over2$}}
\begin{document}

\maketitle

\begin{abstract}
This paper is concerned with the well known Jeffreys--Lindley paradox. In a Bayesian set up, the so-called paradox arises when a point null  hypothesis is tested and an objective prior is sought for the alternative hypothesis. In particular, the posterior for the null hypothesis tends to one when the uncertainty, i.e. the variance, for the parameter value goes to infinity. We argue that the appropriate way to deal with the paradox is to use simple mathematics, and that any philosophical argument is to be regarded as irrelevant.
\end{abstract}

\section{Introduction}\label{sc_intro}
The literature on the \emph{Jeffreys--Lindley paradox} has been prolific since it was brought to the attention of objective Bayesians by \cite{Lindley:1957}. Many authors have discussed this so-called paradox from varying perspectives;  including not only statisticians, but philosophers too. Our aim is to consider the problem using simple mathematics. 

\cite{Lindley:1957} shows that, for point null hypothesis testing, there may be a concern with  the objective Bayesian approach. In the specific example used,  if the prior for the location parameter, in the alternative model to the parameter being zero, has infinite variance, then the Bayesian will always select the null model, regardless of the observed data. This was first suggested as a warning against using improper priors, but the consequences have now become far reaching with a substantial amount of literature written about the observation.

Let us describe the mathematical setting of the problem.
Suppose we wish to test the hypothesis 
$$H_0:\theta=0\quad\mbox{vs}\quad H_1:\theta\ne 0$$ for the normal model $N(x|\theta,1)$. Let $\rho_0=P(M_0)$ be the prior probability assigned to the null hypothesis and let $\pi(\theta)=N(\theta|0,\sigma^2)$, for some $\sigma>0$, be the prior distribution for the unknown parameter $\theta$ under the alternative model. 

Then the Bayes factor for this problem is given by
$$B_{01} = \frac{N(x|0,1)}{\int N(x|\theta,1)\pi(\theta)\;d\theta},$$
which represents the odds in favour of the null hypothesis with respect to the alternative. The decision on whether one  rejects  $H_0$ in favour of  $H_1$ is based on the posterior probability, given by
$$P(M_0|x) = \left[1+\frac{1-P(M_0)}{P(M_0)}\frac{1}{B_{01}}\right]^{-1}.$$
This is the extent of the mathematical foundations to the problem. As Lindley noted, there are some combinations of $(\rho_0,\sigma)$ yielding a $P(M_0|x)$ which one would not wish to countenance. 

The natural objective choice for $\pi(\theta)$ involves taking $\sigma=\infty$. However, rather than a direct plug in of this value, a more general setting has been suggested and considered by \cite{Robert:1993} which is to let $\rho_0$ depend on $\sigma$, i.e. we have $\rho_0(\sigma)$, so then it is possible to study $P(M_0|x)$ as $\sigma\rightarrow\infty$. In this case we can identify three scenarios for $P(M_0|x)$ as $\sigma\rightarrow\infty$, all of which have associated problems. That is, there is no setting, i.e. choice of $\rho_0(\sigma)$, in which the choice $\sigma=\infty$ as an objective choice can work. What we mean by this is explained in Section 2. The conclusion is that the  objective idea of $\sigma=\infty$ does not work and consequently the message is not to use it. On the other hand, we can set the pair $(\sigma<\infty,\rho_0(\sigma))$ objectively using ideas of Type I error calculations and a novel approach to the selection of priors for models. 

The layout of the paper is as follows. In Section \ref{sc_paradox} we formalise the Jeffreys--Lindley paradox and discuss \cite{Robert:1993} solution to it. Section \ref{sc_our} is dedicated to the our approach, and Section \ref{sc_disc} is reserved to conclusions and final comments.

\section{Formalisation of the paradox}\label{sc_paradox}
In order to discuss approaches to the Jeffreys--Lindley paradox, let us first formalise it and, at the same time, define the notation. The objective is to compare the two normal models,
$$M_0=\left\{N(x|0,1) =(2\pi)^{-1/2} \exp(-\half x^2)\right\},$$
$$M_1=\left\{N(x|\theta,1)=(2\pi)^{-1/2}\exp\{-\half(x-\theta)^2\}\right\}.$$
To apply the Bayesian approach, as described in Section \ref{sc_intro}, we need to define both the priors; i.e. the value of $\sigma$, for the unknown parameter $\theta$, and the prior for the null hypothesis; i.e. the value of $\rho_0$. To be most general we will assume that $\rho_0$ can depend on $\sigma$ and hence we write it as $\rho_0(\sigma)$. 

With this information we can compute the Bayes factor representing the odds of the null hypothesis $H_0$. That is
$$B_{01} = \frac{N(x|0,1)}{\int N(x|\theta,1)\cdot N(\theta|0,\sigma^2)\;d\theta}=\frac{e^{-\half x^2}}{e^{-\half x^2/(\sigma^2+1)}}\sqrt{\sigma^2+1},$$
so the posterior probability for the null hypothesis is given by
\begin{equation}\label{eq_intro_2}
P(H_0|x) = \left[1+\frac{1-\rho_0}{\rho_0}\frac{1}{B_{01}}\right]^{-1}=\left[1+\frac{1-\rho_0}{\rho_0}\frac{e^{-\half x^2/(\sigma^2+1)}}{e^{-\half x^2}}\frac{1}{\sqrt{\sigma^2+1}}\right]^{-1}.
\end{equation}

\noindent
We note in (\ref{eq_intro_2}) that the quantity
$$m(\sigma) = \frac{1-\rho_0(\sigma)}{\rho_0(\sigma)}\frac{1}{\sqrt{1+\sigma^2}}$$
is the key term and opens the way to understanding the paradox. We now assume the decision maker wants to select $\sigma=\infty$ in order to implement an objective Bayesian approach. To adequately understand this procedure we argue the decision maker needs to specify $m(\sigma)$ as $\sigma\rightarrow \infty$, and to this end we identify 3 important and exhaustive cases:

\begin{description}

\item (i) $m(\sigma)\rightarrow 0$. 
Under this scenario we have the undesirable result that $P(H_0|x)$ converges to one regardless of the $x$ value. This is the so-called paradoxical result. In fact, $m(\sigma)\rightarrow \infty$
whenever for large $\sigma$ we have, for any $\varepsilon>0$,
\begin{eqnarray*}
\frac{1-\rho_0(\sigma)}{\rho_0(\sigma)}\frac{1}{\sqrt{1+\sigma^2}} &<& \varepsilon \nonumber \\
\mbox{i.e.}\,\,\frac{1-\rho_0(\sigma)}{\rho_0(\sigma)} &\leq& \varepsilon\,\sqrt{1+\sigma^2} \nonumber \\
\rho_0(\sigma) &\geq& \frac{1}{1+\varepsilon\sqrt{1+\sigma^2}}. \nonumber
\end{eqnarray*}
So if the prior on the null hypothesis is too large as $\sigma\rightarrow \infty$; i.e. $\sigma\rho_0(\sigma)\rightarrow \infty$,  then the posterior probability on the null hypothesis will converge to 1. 

\item (ii) $m(\sigma)\rightarrow c$ for some constant $0<c<\infty$. Under this scenario
it is that, for large $\sigma$,
$$\rho_0(\sigma)=\frac{1}{1+c\sqrt{1+\sigma^2}}\approx\frac{1}{1+c\sigma}.$$ 
In particular, \cite{Robert:1993} presents an objective argument for 
$$\rho_0(\sigma) = \frac{1}{1+\sqrt{2\pi}\sigma}.$$
However, this idea leads to an undesirable inconsistency in that $\rho_0(\sigma)\rightarrow 0$ yet $P(M_0|x)$ is converging to a constant bounded away from 0. Thus, with $\sigma=\infty$, we have $P(M_0)=0$ but $P(M_0|x)\ne 0$, which are incoherent choices.
 
\item (iii) $m(\sigma)\rightarrow\infty$. Under this scenario we have that $P(M_0|x)\rightarrow 0$. This at least now becomes consistent with the prior probability since $\rho_0(\sigma)\rightarrow 0$ in this case. Yet undesirable in that with $\sigma=\infty$, $P(M_0|x)=0$.

\end{description}

These considerations clearly exclude the choice $\sigma=\infty$. It simply does not work. Thus a finite choice of $\sigma$ is required. In the next section we will demonstrate how we can set $(\sigma<\infty,\rho_0(\sigma))$ objectively.  

\section{An objective choice for ($\sigma,\rho_0(\sigma)$)}\label{sc_our}
Given a value of $\sigma$ we first, in Section \ref{sc_rho}, show how to obtain an objective choice for $\rho_0(\sigma)$. Then, in Section \ref{sc_sigma}, we show how $\sigma<\infty$ can be selected objectively.

\subsection{The prior $\rho_0(\sigma)$}\label{sc_rho}
Our approach consists in measuring the \emph{worth} of the alternative hypothesis with respect to the null, as outlined in \cite{Villa:Walker:2014}. In particular, we apply the well known asymptotic Bayesian property that, if a model is misspecified, the posterior accumulates at the model which the nearest, in terms of Kullback--Leibler divergence, to the true model \citep{Berk:1966}. As such, the divergence $D_{KL}(N(x|\theta,1)\|N(x|0,1))$ represents the loss we would incur if model $M_1$ is removed and it is true. Since we do not know $\theta$, but we have the prior $\pi(\theta)$, we can compute the expected loss as
\begin{equation}\label{eq_ourprior1}
\int_{\Theta}D_{KL}\bigg(N(x|\theta,1)\|N(x|0,1)\bigg) \pi(\theta) \; \mbox{d}\theta = \int \half \theta^2\,\pi(\theta)\;\mbox{d}\theta=\half\sigma^2.
\end{equation}
The model prior is determined by means of the \emph{self-information} loss function \citep{Merhav:Feder:1998}, which represents the loss connected to a probability statement. For model $M$, the self-information loss is given by $-\log P(M)$. Therefore, by equating the self-information with the expected loss determined in (\ref{eq_ourprior1}), we have that the prior on the alternative model is
$$1-\rho_0(\sigma) \propto e^{\sigma^2}.$$
Note that the prior for the null hypothesis is $\rho_0(\sigma)\propto1$, and so we have
$$\rho_0(\sigma)=\frac{1}{1+\exp\{\half\sigma^2\}}.$$
This then fits into category (iii) for large $\sigma$, which implies that $P(M_0|x,\sigma)$ goes to zero  as $P(M_0)\rightarrow 0$. Thus there is coherence in this approach; however, we are not advocating the choice of large $\sigma$.

\subsection{Determining $\sigma$}\label{sc_sigma}
In any classical test the Type I error is of key importance. We can use this quantity to objectively set the value for $\sigma$; if indeed the Type I error is an objective quantity, but nevertheless it needs to be set, and a valid objective Bayesian criterion is to match classical benchmarks and quantities. 

To determine an appropriate value for $\sigma$ based on the classical concept of Type I error, we would select $\sigma$ so that
$$P_0(\mbox{reject } H_0) = \alpha,$$
where $\alpha\in(0,1)$ and $P_0$ is the probability under the null hypothesis. Regardless of the surroundings, all Bayesian experimenters in this problem would need to assign an $\alpha_B$ value for which one would reject $H_0$ if $P(M_0|x)<\alpha_B$. 
To have
\begin{equation}
P_0\bigg(P(M_0|x)<\alpha_B\bigg) = \alpha,\label{nedo}
\end{equation}
we require
\begin{eqnarray}\label{eq_sigma_1}
\frac{1}{1+m(\sigma)\exp\left\{\half x^2\frac{\sigma^2}{1+\sigma^2}\right\}} &<& \alpha_B, \nonumber \\
\mbox{i.e.}\,\,\exp\left\{\frac{1}{2} x^2\frac{\sigma^2}{1+\sigma^2}\right\} &>& \frac{1/\alpha_B-1}{m(\sigma)} \nonumber \\
\frac{1}{2}x^2\frac{\sigma^2}{1+\sigma^2} &>& \log\left(\frac{\alpha_B^{-1}-1}{m(\sigma)}\right) \nonumber \\
x^2 &>& \frac{2(1+\sigma^2)}{\sigma^2}\log\left(\frac{\alpha_B^{-1}-1}{m(\sigma)}\right).
\end{eqnarray}
Therefore, if we write
$$\psi(\sigma) = \frac{2(1+\sigma^2)}{\sigma^2}\log\left(\frac{\alpha_B^{-1}-1}{m(\sigma)}\right),$$
we have (\ref{nedo}) as
\begin{equation}\label{eq_sigma_2}
P_0\bigg(x^2>\psi(\sigma)\bigg) = 2\left[1-\Phi\left(\sqrt{\psi(\sigma)}\right)\right]=\alpha.
\end{equation}

\begin{figure}[]
\centering
\includegraphics[scale=0.5]{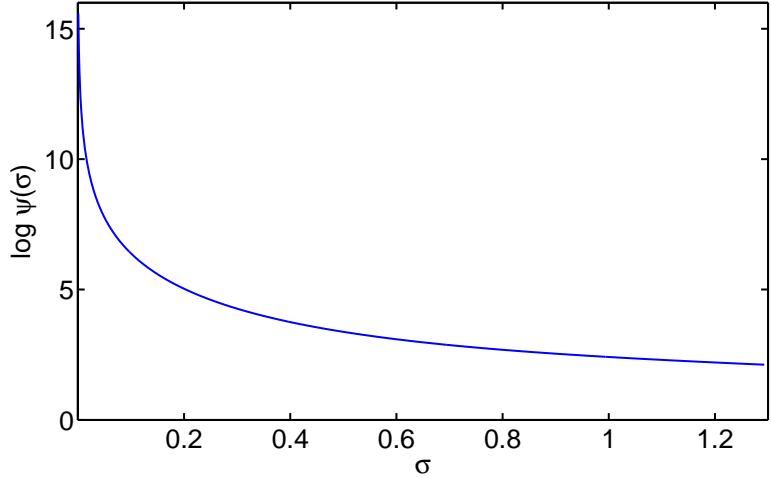}
\caption{Plot of $\log \psi(\sigma)$, with $\alpha_B=0.05$, for $\sigma<1.3$. By setting $\sigma=1.2933$ we ensure that $\psi(\sigma)>0$.}
\label{fig1}
\end{figure}

The key here is that $\psi(\sigma)$ is decreasing as $\sigma$ increases, so there is a one-to-one correspondence between $\alpha$ and $\sigma$ satisfying (\ref{eq_sigma_2}). Figure \ref{fig1} shows the behaviour of $\log \psi(\sigma)$, given $\alpha_B=0.05$. As it must be that  $\psi(\sigma)>0$, we compute $\log\psi(\sigma)$ up to $\sigma=1.2930$, which is the value that ensures $m(\sigma)<\alpha_B^{-1}-1$, therefore, a positive $\psi(\sigma)$.

Expression (\ref{eq_sigma_2}) has to be solved numerically. So, for example, if $\alpha_B=\alpha=0.05$, we would have $\sigma=0.44$. In other words, we can be objective about $\sigma$ with a finite value. The notion therefore that an objective $\sigma$ and $\sigma=\infty$ is the only choice is wrong.  An objective classical test requires an $\alpha$ value and it is this which can be linked to the (finite) objective choice for $\sigma$. 

\section{Discussion}\label{sc_disc}
The findings of this paper can be summarised as follows. The posterior for the point null hypothesis is driven by the quantity $m(\sigma)$; in particular, if $\sigma=\infty$ is desired as an objective criterion then the  behaviour of $m(\sigma)$ as $\sigma\rightarrow\infty$ is the key. If the prior $\rho_0(\sigma)$ is fixed, e.g. is equal to $\half$, then the Jeffreys--Lindley paradox arises, since the posterior probability $P(H_0|x)$ goes to one. \cite{Robert:1993} proposed to solve the issue by having $m(\sigma)$ to converge to a positive constant. Although the direct paradox is avoided, the approach gives an incoherent result as the posterior mass on $H_0$ is positive whereas the prior mass is zero. Our approach gives a quantity $m(\sigma)$ which goes to infinity, for $\sigma$ going to infinity, which both solves the paradox and yields zero posterior mass for $H_0$ when $\rho_0(\sigma)=0$, implying the prior mass for $H_0$ is zero.

It is clear that the three types of behaviour of $m(\sigma)$ for large $\sigma$ rule out the possibility of having $\sigma=\infty$. As such, $\sigma$ has to be determined to have a finite value. For $\rho_0(\sigma)$, the choice can be either objective or subjective. Our approach allows  $\rho_0(\sigma)$ to be determined in an objective fashion by considering the loss in information if the true model is removed. \cite{Dellap:2012}, on the other hand, propose a prior for the null hypothesis that is subjective.

\cite{Dellap:2012} focus on models for which the use of a multivariate normal prior is appropriate, such as linear regression models, generalised linear models and standard time series models. The idea is to set the multiplicative constant for the prior dispersion matrix, $c_m$, which will indicate the level of prior uncertainty. The authors aim to reduce the sensitivity of the posterior model probabilities to the scale of the prior by suitably specifying the prior model probability. This is done by setting
$$P(M) \propto P^\prime(M)c_M^{d_m},$$
where $d_M$ is the dimension of the model $M$ and $P^\prime(M)$ is a suitably determined base line prior model probability. \cite{Dellap:2012} recommend $P^\prime(M)\propto1$, although other choices are possible. We see that the core of the whole approach is to make a prior model probability dependent on the variance of the prior in the parameters, avoiding the Jeffreys--Lindley paradox. \\

The conclusion is that it is not possible to be objective for $\pi(\theta)$ by setting $\sigma=\infty$. This is not the sole case where objective Bayes fails to deliver adoptable solutions. For example, Jeffreys' rule prior for multidimensional parameter spaces gives prior distribution with poor performance properties \citep{BerSmi:1994}. It is common practice not to use Jeffreys prior in these type of problems and opt for a different solution, such as reference priors.

However, an objective and finite value of $\sigma$ can be assigned by exploiting thinking behind classical tests and setting the Type I error. That is, there is a one-to-one correspondence between $\sigma$ and the Type I error $\alpha$ and it is this correspondence which permits the  interpretation and assignment of $\sigma$.  

Surprisingly, or not, there have been philosophical papers attempting to find some hidden profound explanation behind the paradox; see, for example, the recent papers of \cite{Spanos:2013} and \cite{Sprenger:2013}. We argue that it is not necessary to philosophize, as the mathematics of the problem are quite straightforward and a clear picture of what is happening can be understood solely by mathematical considerations. 

To discuss some of the philosophy, \cite{Spanos:2013} says: ``The question that generally arises is why the Bayesian and the likelihoodist approaches give rise to the above conflicting and confusing results". However, we have $P(M_0|x)<\alpha_B \Leftrightarrow x^2>\psi(\sigma)$, which is precisely the form of the classical test !  

The classical test is: reject $H_0$ if $x^2>c_\alpha$, where $P_0(x^2>c_\alpha)=\alpha$. We can then set $\psi(\sigma)=c_\alpha$ to ensure a standard value for the Type I error. Consequently, Bayes makes no contribution to this problem, since even a subjective Bayesian approach will yield a classical test, but with perhaps a non-standard Type I error. Such an observation between Bayesian and classical tests has been made by \cite{ShiWalker:2013}.

In short, both Bayesian and classical tests reject $H_0$ if $x^2>c$, and this is the obvious procedure for testing $H_0:\theta=0$. How one determines $c$ makes the difference, either via a Type I error, $\alpha$, or via a prior $\pi(\theta)$, i.e. $\sigma$, but nevertheless there is a one-to-one correspondence between the two.

\end{document}